\documentclass[10pt,a4paper]{article}
\usepackage{amsfonts}
\usepackage{}
\usepackage{amsmath}
\usepackage{latexsym}
\usepackage{amssymb}
\usepackage{amstext}
\usepackage{fancyhdr}
\usepackage{indentfirst}
\usepackage{graphicx}
\usepackage{epsfig}
\usepackage{pgf}
\usepackage{listings}
\usepackage{cases}
\usepackage{epstopdf}
\usepackage{graphicx}
\usepackage{enumerate}
\date{}
\makeatletter

\newcommand{\Rmnum}[1]{\expandafter\@slowromancap\romannumeral #1@}
\makeatother
\usepackage{graphicx}
\usepackage{float}
\usepackage{subfigure}
\usepackage{bm}
\begin{document}
\title{More on minors of Hermitian (quasi-)Laplacian matrix of the second kind for mixed graphs}
\author{Qi Xiong$^a$, Gui-Xian Tian$^{a}$\footnote{Corresponding author. E-mail addresses: 252321425@qq.com (Q. Xiong), gxtian@zjnu.cn (G.-X. Tian), cuishuyu@zjnu.cn (S.-Y. Cui).}, Shu-Yu Cui$^b$\\\\%EndAName
    {\small{\it $^a$Department of Mathematics,}}
    {\small{\it Zhejiang Normal University, Jinhua, 321004, China}}\\
    {\small{\it $^b$Xingzhi College, Zhejiang Normal University, Jinhua, 321004, China}}
}\maketitle

\begin{abstract} A mixed graph $M_{G}$ is the graph obtained from an unoriented simple graph $G$ by giving directions to some edges of $G$, where $G$ is often called the underlying graph of $M_{G}$. In this paper, we introduce two classes of incidence matrices of the second kind of $M_{G}$, and discuss the determinants of these two matrices for rootless mixed trees and unicyclic mixed graphs. Applying these results, we characterize the explicit expressions of various minors for Hermitian (quasi-)Laplacian matrix of the second kind of $M_{G}$. Moreover, we give two sufficient conditions that the absolute values of all the cofactors of Hermitian (quasi-)Laplacian matrix of the second kind are equal to the number of spanning trees of the underlying graph $G$.


\end{abstract}

\section{Introduction}

The classical matrix tree theorem, attributed to Kirchhoff \cite{Kirchhoff}, states that the number of spanning trees of an unoriented graph can be acquired directly by computing the cofactors of its Laplacian matrix. An elementary proof of the matrix tree theorem can be found in \cite{Cvetkovic}. That is a natural generalization of the matrix tree theorem to oriented graphs, see \cite{De Leenheer} and the cited references therein. Since the oriented graph is a special case of a mixed graph. Then the generalized matrix tree theorem of mixed graphs was investigated in terms of a real Laplacian matrix of mixed graphs in \cite{Bapat}. Recently, the concept of Hermitian adjacency matrix of a mixed graph was first proposed by Guo and Mohar \cite{Guo}, and Liu and Li \cite{Liu} independently. And then Hermitian (quasi-)Laplacian matrix of mixed graphs was introduced in \cite{Yu1,Yu2}. Principal minor version of matrix tree theorem for mixed graphs was surveyed in \cite{Sarma,Yu3}. In 2020, Hermitian adjacency matrix of the second kind for mixed graphs was introduced by Mohar \cite{Mohar}, and more properties about this matrix can be found in \cite{Li}. This paper mainly focuses on minors of Hermitian (quasi-)Laplacian matrix of the second kind for mixed graphs. We characterize the expressions of various minors for Hermitian (quasi-)Laplacian matrix of the second kind. Our main strategy is to decompose the Hermitian (quasi-)Laplacian matrix of the second kind into a product of an incidence matrix and its conjugate transpose, and use the Cauchy-Binet Theorem to calculate various minors for Hermitian (quasi-)Laplacian matrix of the second kind. Our strategy is not a novel idea and has been used in several references, see \cite{Bapat,Cvetkovic,Sarma}. However, we find that, for some special types of mixed graphs, the number of spanning trees of their underlying graphs can be obtained by computing the cofactors of Hermitian (quasi-)Laplacian matrix of the second kind. Moreover, we give two sufficient conditions that the absolute values of all the cofactors of Hermitian (quasi-)Laplacian matrix of the second kind for a mixed graph are equal to the number of spanning trees of its underlying graph.

All the graphs discussed are simple and finite connected graphs throughout this paper. Let $G=(V,E)$ be a graph with vertex set $V=\{v_1,v_2,\ldots, v_n\}$ and edge set $E=\{e_1,e_2,\ldots,e_m\}$. We call $v_i$ is adjacent to $v_j$ in $G$ if there is an edge $e$ connecting $v_i$ and $v_j$. A mixed graph $M_{G}$ is obtained from $G$ by giving directions to some edges of $G$, where $G$ is often called the underlying graph of $M_{G}$. If the direction of $e$ is from $v_i$ to $v_j$ in $M_{G}$, then we write it as $v_i\rightarrow v_j$. If the direction is reversed, then we write it as $v_i\leftarrow v_j$, or $v_j\rightarrow v_i$. If $e$ is an undirected edge in $M_{G}$, then we write it as $u\leftrightarrow v$.

Given a mixed graph $M_{G}$, the Hermitian adjacency matrix of the second kind for $M_{G}$ is the $n\times n$ matrix $N(M_{G})=(n_{ij})$, whose entries are given by
\begin{equation*}
 n_{ij}=\left\{\begin{array}{rcl}
&1,&\quad\text{if}\; v_i\leftrightarrow v_j;\\
&\frac{1+\sqrt{3}\mathbf{i}}{2},&\quad\text{if}\; v_i\rightarrow v_j;\\
&\frac{1-\sqrt{3}\mathbf{i}}{2},&\quad\text{if}\; v_i\leftarrow v_j;\\
&0,&\quad\text{otherwise}.
\end{array}\right.
\end{equation*}
Let $D(M_{G})$ be the degree diagonal matrix of the underlying graph $G$ of $M_{G}$. Then we introduce the Hermitian Laplacian matrix and Hermitian quasi-Laplacian matrix of the second kind for $M_{G}$, which are defined as $L(M_{G})=D(M_{G})-N(M_{G})$ and $Q(M_{G})=D(M_{G})+N(M_{G})$, respectively.

It is well known that, for any unoriented graph, the absolute value of any cofactor of its Laplacian matrix is indeed equal to the number of its spanning trees. Our motivation is the question that for which mixed graphs this property still holds related to the Hermitian (quasi-)Laplacian matrix of the second kind.

\begin{figure}[H]
\centering
\subfigure[$G$]{
\begin{minipage}[t]{0.32\linewidth}
\centering
\includegraphics[width=1.1in]{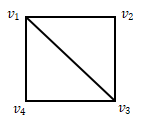}
\end{minipage}}
\subfigure[$M_{G}'$]{
\begin{minipage}[t]{0.32\linewidth}
\centering
\includegraphics[width=1.1in]{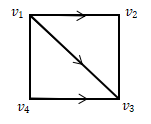}
\end{minipage}}
\subfigure[$M_{G}''$]{
\begin{minipage}[t]{0.32\linewidth}
\centering
\includegraphics[width=1.1in]{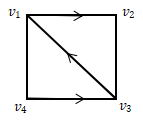}
\end{minipage}}
\centering
\caption{Two mixed graphs $M_{G}'$ and $M_{G}''$ with the underlying graph $G$ in Example 1.1.}
\end{figure}
\noindent
\textbf{Example 1.1.} Let $G$ be the unoriented graph with 4 vertices as shown in Figure 1(a). Then the Laplacian matrix of $G$ is given by
\begin{equation*}
L(G)=\left(\begin{array}{cccc}
3 &-1&-1&-1\\
-1&2&-1&0\\
-1&-1&3&-1\\
-1&0&-1&2\\
\end{array}\right).
\end{equation*}
Let $L_{i,j}$ be a submatrix after deleting row $i$ and column $j$ from $L(G)$. It is easy to verify that $|\det L_{i,j}(G)|=8$ for $1\leq i,j\leq 4$. So the number of spanning trees of $G$ is equal to $8$.

For the mixed graph $M_{G}'$ as shown in Figure 1(b), its Hermitian Laplacian matrix of the second kind is given by
\begin{equation*}
L(M_{G}')=\left(\begin{array}{cccc}
3 &-\frac{1+\sqrt{3}\mathbf{i}}{2}&-\frac{1+\sqrt{3}\mathbf{i}}{2}&-1\\
-\frac{1-\sqrt{3}\mathbf{i}}{2}&2&-1&0\\
-\frac{1-\sqrt{3}\mathbf{i}}{2}&-1&3&-\frac{1-\sqrt{3}\mathbf{i}}{2}\\
-1&0&-\frac{1+\sqrt{3}\mathbf{i}}{2}&2\\
\end{array}\right).
\end{equation*}
Then by calculation we get $|\det L_{i,j}(M_{G}')|=8$ for $1\leq i,j\leq 4$, which implies that the number of spanning trees of the underlying graph $G$ is equal to the absolute value of any cofactor of the Hermitian Laplacian matrix of the second kind for $M_{G}'$. However, for the mixed graph $M_{G}''$ as shown in Figure 1(c), its Hermitian Laplacian matrix of the second kind is given by
\begin{equation*}
L(M_{G}'')=\left(\begin{array}{cccc}
3 &-\frac{1+\sqrt{3}\mathbf{i}}{2}&-\frac{1-\sqrt{3}\mathbf{i}}{2}&-1\\
-\frac{1-\sqrt{3}\mathbf{i}}{2}&2&-1&0\\
-\frac{1+\sqrt{3}\mathbf{i}}{2}&-1&3&-\frac{1-\sqrt{3}\mathbf{i}}{2}\\
-1&0&-\frac{1+\sqrt{3}\mathbf{i}}{2}&2\\
\end{array}\right).
\end{equation*}
It is easy to check that the modules of different cofactors are not equal, and they are not equal to the number of spanning trees of the underlying graph $G$.\\

Note that for unoriented graph, there is generally no such property that the absolute value of any cofactor of its quasi-Laplacian matrix is equal to the number of its spanning trees. But this property holds for the Hermitian quasi-Laplacian matrix of the second kind of some special mixed graphs.

\begin{figure}[H]
\centering
\subfigure[$G$]{
\begin{minipage}[t]{0.32\linewidth}
\centering
\includegraphics[width=1.1in]{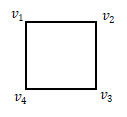}
\end{minipage}}
\subfigure[$M_{G}'$]{
\begin{minipage}[t]{0.32\linewidth}
\centering
\includegraphics[width=1.1in]{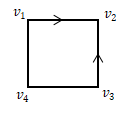}
\end{minipage}}
\subfigure[$M_{G}''$]{
\begin{minipage}[t]{0.32\linewidth}
\centering
\includegraphics[width=1.1in]{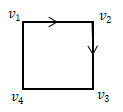}
\end{minipage}}
\centering
\caption{Two mixed graphs $M_{G}'$ and $M_{G}''$ with the underlying graph $G$ in Example 1.2.}
\end{figure}

\textbf{Example 1.2.} Let $G$ be an unoriented graph with 4 vertices as shown in Figure 2(a). Also let $M_{G}'$ and $M_{G}''$ be two mixed graphs obtained by orienting some edges of $G$ as shown in Figure 2(b) and Figure 2(c). It is clear that the number of spanning trees of $G$ equals 4.

For the mixed graph $M_{G}'$, we get its Hermitian quasi-Laplacian matrix of the second kind in the following.
\begin{equation*}
Q(M_{G}')=\left(
\begin{array}{cccc}
2 &\frac{1+\sqrt{3}\mathbf{i}}{2}&0&1\\
\frac{1-\sqrt{3}\mathbf{i}}{2}&2&\frac{1-\sqrt{3}\mathbf{i}}{2}&0\\
0&\frac{1+\sqrt{3}\mathbf{i}}{2}&2&1\\
1&0&1&2\\
\end{array}
\right).
\end{equation*}
Then by calculation we can get that $|\det Q_{i,j}(M_{G}')|=4$ for $1\leq i,j\leq 4$, which implies that the absolute value of any cofactor of the Hermitian quasi-Laplacian matrix of the second kind for $M_{G}'$ is equal to the number of spanning trees of the underlying graph $G$. However for the mixed graph $M_{G}''$, its Hermitian quasi-Laplacian matrix of the second kind can be written as follows.
\begin{equation*}
Q(M_{G}'')=\left(
\begin{array}{cccc}
2 &\frac{1+\sqrt{3}\mathbf{i}}{2}&0&1\\
\frac{1-\sqrt{3}\mathbf{i}}{2}&2&\frac{1+\sqrt{3}\mathbf{i}}{2}&0\\
0&\frac{1-\sqrt{3}\mathbf{i}}{2}&2&1\\
1&0&1&2\\
\end{array}
\right).
\end{equation*}
Then by calculation we find that the modules of different cofactors are not equal. This implies that the number of spanning trees of the underlying graph $G$ can not be calculated in this way.

The two examples above show that the matrix tree theorem still holds for some special mixed graphs in terms of the Hermitian (quasi-)Laplacian matrix of the second kind. Therefore, we shall characterize the mixed graphs with this property in Section 3. In order to achieve our goals, we first give two classes of incidence matrices of the second kind of mixed graphs. Then the determinants of these two matrices are discussed for rootless mixed trees and unicyclic mixed graphs. Applying these results, we give expressions of principal minors of the Hermitian (quasi-)Laplacian matrix of the second kind for a mixed graph. These contents are included in Section 2. Section 3 deals with non-principal minors for the Hermitian (quasi-)Laplacian matrix of the second kind of mixed graphs. Finally, we give two sufficient conditions that the absolute values of all the cofactors of Hermitian (quasi-)Laplacian matrix of the second kind are equal to the number of spanning trees of the underlying graph $G$.

\section{Principal minors}
Recall that for a mixed graph $M_{G}$, its incidence matrix of the second kind, denoted by $S(M_{G})=(s_{ie})$, was introduced in \cite{Xiong}. It is defined as
\begin{equation*}
s_{ie}=\left\{\begin{array}{rcl}
&-s_{je},&\quad\text{if}\; v_{i}\leftrightarrow v_{j};\\
&-\frac{1+\sqrt{3}\textbf{i}}{2}s_{je},&\quad \text{if}\; v_{i}\rightarrow v_{j};\\
&-\frac{1-\sqrt{3}\textbf{i}}{2}s_{je},&\quad \text{if}\; v_{i}\leftarrow v_{j};\\
&0,&\quad\text{otherwise,}
\end{array} \right.
\end{equation*}
where $s_{ie}$ is a complex number satisfying $|s_{ie}|=1$ or 0. Notice that $S(M_{G})$ is not unique. For simplicity sake, we provide a special case of this matrix. We assign every undirected edge of $M_{G}$ to any direction and denote this new graph by $M_{G'}$. Then the Hermitian incidence matrix of the second kind $S(M_{G})=(s_{ie})$ can be rewritten as
\begin{equation*}
 s_{ie}=\left\{
\begin{array}{rcl}
&1,& \text{if  $e$ is a new directed edge of $M_{G'}$ with $v_{i}$ as its head};\\
&-1,& \text{if $e$ is a new directed edge of $M_{G'}$ with $v_{i}$ as its tail};\\
&1,& \text{if $e$ is a directed edge of $M_{G}$ with $v_{i}$ as its head};\\
&-\frac{1-\sqrt{3}\mathbf{i}}{2},& \text{if $e$ is a directed edge of $M_{G}$ with $v_{i}$ as its tail};\\
&0,& \text{otherwise}.
\end{array} \right.
\end{equation*}
Similarly, the quasi-incidence matrix of the second kind $T(M_{G})=(t_{ie})$ \cite{Xiong} is given by
\begin{equation*}
 t_{ie}=\left\{\begin{array}{rcl}
&t_{je},&\text{if}\quad v_{i}\leftrightarrow v_{j};\\
&\frac{1+\sqrt{3}\mathbf{i}}{2}t_{je},& \text{if}\quad v_{i}\rightarrow v_{j};\\
&\frac{1-\sqrt{3}\mathbf{i}}{2}t_{je},& \text{if}\quad v_{i}\leftarrow v_{j};\\
&0,&\text{otherwise},
\end{array} \right.
\end{equation*}
where $t_{ie}$ is a complex number satisfying $|t_{ie}|=1$ or 0. Since this matrix is also non-unique, then we also give a special case of this incidence matrix as a matter of convenience. Let $T(M_{G})=(t_{ie})$, in which the elements in this matrix are as follows.
\begin{equation*}
 t_{ie}=\left\{\begin{array}{rcl}
&1,& \text{if $e$ is an undirected edge of $M_{G}$ with $v_{i}$ as its end vertex};\\
&1,&\text{if $e$ is a directed edge of $M_{G}$ with $v_{i}$ as its head};\\
&\frac{1-\sqrt{3}\mathbf{i}}{2},& \text{if $e$ is a directed edge of $M_{G}$ with $v_{i}$ as its tail};\\
&0,&\text{otherwise}.
\end{array} \right.
\end{equation*}
It was proved \cite{Xiong} that $L(M_{G})=S(M_{G})S(M_{G})^{\ast}$ and $Q(M_{G})=T(M_{G})T(M_{G})^{\ast}$ are both positive semidefinite matrices, where $A^{\ast}$ denotes the conjugate transpose of a matrix $A$.

Let $M_{G}$ be a mixed graph with vertex set $V$ and edge set $E$, and $V_{1}\subseteq V$, $E_{1}\subseteq E$, then $M_{R}$ is a substructure of $M_{G}$ with vertex subset $V_{1}$ and edge subset $E_{1}$. We call $M_{R}$ a square substructure if the number of vertices and edges is the same. A rootless tree is a mixed tree obtained by deleting only one vertex in the tree, and this vertex is called the root of rootless tree. Let $M_{T_{u}}$ be a rootless tree with root vertex $u$, and a directed edge $e=vw$ of $M_{T_{u}}$ with the direction from $v$ to $w$. We call $e$ is away (resp. towards) the root vertex $u$ if $d_{uv}=d_{uw}-1$ (resp. $d_{uv}=d_{uw}+1$) in the tree $M_{T_{u}}$, where $d_{uw}$ denotes the distance between two vertices $u$ and $w$ in the underlying tree of $M_{T_{u}}$.\\
\\
\textbf{Lemma 2.1.}  Suppose that $M_{G}$ is a mixed graph and a substructure $M_{R}$ of $M_{G}$ is a rootless mixed tree. Then
\[
\det(S(M_{R}))=(-\frac{1-\sqrt{3}\mathbf{i}}{2})^{\alpha}(-1)^{\beta},
\]
where $\alpha$ is the number of the old directed edges away from the root in $M_{R}$, and $\beta$ is the number of the new directed edges away from the root after we assign any direction to all the undirected edges in $M_{R}$.\\
\\
\textit{Proof.} First suppose that the substructure $M_{R}$ is a rootless mixed path with vertex set $\{v_{1},v_{2},\ldots,v_{n}\}$ and edge set $\{e_{1},e_{2},\ldots,e_{n}\}$. After a suitable relabeling of the vertices if necessary, we may assume that the edge $e_{n}$ has only one endpoint $v_{n}$ and $e_{k}=v_{k}v_{k+1}$ for $k=1,2,\ldots,n-1$. Thus, the incidence matrix of the second kind of $S(M_{R})$ can be written as
\begin{equation*}
S(M_{R})=\left(
\begin{array}{ccccccc}
a_1 &0&0&\cdots&0&0&0\\
b_1 &a_2&0&\cdots&0&0&0\\
0 &b_2&a_3&\cdots&0&0&0\\
\vdots&\vdots&\vdots&\ddots&\vdots&\vdots&\vdots\\
0&0&0&\cdots&a_{n-2}&0&0\\
0&0&0&\cdots&b_{n-2}&a_{n-1}&0\\
0&0&0&\cdots&0&b_{n-1}&a_n\\
\end{array}
\right).
\end{equation*}

According to the definition of incidence matrix of the second kind, we obtain that $a_i$ could be $1$, $-1$, or $-\frac{1-\sqrt{3}\mathbf{i}}{2}$. Furthermore, the numbers of $-\frac{1-\sqrt{3}\mathbf{i}}{2}$ and $-1$ on the diagonal are equal to the number $\alpha$ of old directed edges away from the root and the number $\beta$ of the new directed edges away from the root after we assign any direction to all the undirected edges in $M_{R}$, respectively. Hence the result follows.

In the following, we assume that the substructure $M_{R}$ is a rootless mixed tree with a hanging edge $e_{n}=v_{n-1}v_{n}$ with one degree vertex $v_n$. Let $S'$ denote the principal submatrix obtained by deleting the row of $S(M_{R})$ corresponding to $v_{n}$ and the column of $S(M_{R})$ corresponding to $e_{n}$. Then $\det(S(M_{R}))=a_n\cdot \det(S')$, in which $a_n$ could be $1$, $-1$, or $-\frac{1-\sqrt{3}\mathbf{i}}{2}$. Repeat the above proceeding to delete all hanging edges until there are only some rootless paths left, implying that the result holds. $\Box$\\
\\
\textbf{Lemma 2.2.} Let $M_{G}$ be a mixed graph with the substructure $M_{R}$. If $M_{R}$ is a rootless mixed tree, then $\det(T(M_{R}))=(\frac{1-\sqrt{3}\mathbf{i}}{2})^{\alpha}$, in which $\alpha$ is the number of directed edges away from the root in $M_{R}$.\\
\\
\textit{Proof.} The proof is exactly similar to that of Lemma 2.1, the details are ignored. $\Box$\\

Given a mixed cycle $M_{C}$ with vertex set $\{v_{1},v_{2},\ldots,v_{n}\}$, in which these vertices are arranged in the clockwise direction. We use $a(M_{C})$, $b(M_{C})$, $c(M_{C})$ to represent the number of directed edges in the clockwise direction, directed edges in the counterclockwise direction and undirected edges in $M_{C}$, respectively. If there is no ambiguity, we mark them $a,b,c$ for the convenience of narration. Then we divide all mixed cycles into the following four categories $\{\Phi_{1},\Phi_{2},\Phi_{3},\Phi_{4}\}$. If a mixed cycle satisfies $a-b\equiv 1$ or $5 (\mod 6)$, then it belongs to $\Phi_{1}$. If a mixed cycle satisfies $a-b\equiv 2$ or $4(\mod 6)$, then it belongs to $\Phi_{2}$. If a mixed cycle satisfies $a-b\equiv 3(\mod 6)$, then it belongs to $\Phi_{3}$. If a mixed cycle satisfies $a-b\equiv 0(\mod 6)$, then it belongs to $\Phi_{4}$. From the definition given above, we can get the following results about mixed cycles.\\
\\
\textbf{Lemma 2.3.} For a mixed cycle $M_{C}$ with $n$ vertices, we have
\begin{equation*}
|\det(S(M_{C}))|=\left\{\begin{array}{rcl}
&1,&\text{if $M_{C}$ belongs to $\Phi_{1}$};\\
&\sqrt{3},&\text{if $M_{C}$ belongs to $\Phi_{2}$};\\
&2,&\text{if $M_{C}$ belongs to $\Phi_{3}$};\\
&0,&\text{if $M_{C}$ belongs to $\Phi_{4}$}.
\end{array} \right.
\end{equation*}
\\
\textit{Proof.} According to the definition of incidence matrix of the second kind, we need to assign any undirected edge of $M_{C}$ to a direction. However, no matter how we orient the undirected edges, the absolute value of the determinant of the matrix $S(M_{C})$ is not going to change. Thus, we might as well give each undirected edge a forward direction. With a proper label of the vertices in $M_{C}$, the matrix $S(M_{C})$ can be written as
\begin{equation*}
S(M_{C})=\left(
\begin{array}{ccccccc}
a_1 &0&0&\cdots&0&0&b_1\\
b_2 &a_2&0&\cdots&0&0&0\\
0 &b_3&a_3&\cdots&0&0&0\\
\vdots&\vdots&\vdots&\ddots&\vdots&\vdots&\vdots\\
0&0&0&\cdots&a_{n-2}&0&0\\
0&0&0&\cdots&b_{n-1}&a_{n-1}&0\\
0&0&0&\cdots&0&b_n&a_n\\
\end{array}
\right).
\end{equation*}
Let $a,b,c$ be as described in the previous paragraph of Lemma 2.3. Then expanding the determinant of $S(M_{C})$ along the top row, one gets that
\begin{equation*}
\begin{aligned}
|\det(S(M_{C}))|=&|(-\frac{1-\sqrt{3}\mathbf{i}}{2})^{b}+(-1)^{n-1}(-1)^{c}(-\frac{1-\sqrt{3}\mathbf{i}}{2})^{a}|\\
=&|(-\frac{1-\sqrt{3}\mathbf{i}}{2})^{b-a}+(-1)^{n+c-1}|\\
=&|(-1)^{b+a}||(\frac{1-\sqrt{3}\mathbf{i}}{2})^{b-a}+(-1)^{n+c-1-b-a}|\\
=&|(\frac{1-\sqrt{3}\mathbf{i}}{2})^{b-a}-1|.\\
\end{aligned}
\end{equation*}
First of all if $M_{C}$ belongs to $\Phi_{1}$, then $a-b\equiv 1$ or $5(\mod 6)$. Thus one obtains that
\begin{equation*}
\begin{aligned}
|\det(S(M_{C}))|=|(\frac{1-\sqrt{3}\mathbf{i}}{2})^{b-a}-1|=|(\frac{1\pm\sqrt{3}\mathbf{i}}{2})-1|=|\frac{-1\pm\sqrt{3}\mathbf{i}}{2}|=1.\\
\end{aligned}
\end{equation*}
Now if $M_{C}$ belongs to $\Phi_{2}$, in this case, $a-b\equiv 2$ or $4(\mod 6)$, then
\begin{equation*}
\begin{aligned}
|\det(S(M_{C}))|=|(\frac{1-\sqrt{3}\mathbf{i}}{2})^{b-a}-1|=|(\frac{-1\pm\sqrt{3}\mathbf{i}}{2})-1|=|\frac{-3\pm\sqrt{3}\mathbf{i}}{2}|=\sqrt{3}.\\
\end{aligned}
\end{equation*}
Again if $M_{C}$ belongs to $\Phi_{3}$, then $a-b\equiv 3(\mod 6)$. So we get that
\begin{equation*}
\begin{aligned}
&|\det(S(M_{C}))|=|(\frac{1-\sqrt{3}\mathbf{i}}{2})^{b-a}-1|=|-1-1|=2.\\
\end{aligned}
\end{equation*}
Finally if $M_{C}$ belongs to $\Phi_{4}$, then $a-b\equiv 0(\mod 6)$. It follows that
\begin{equation*}
\begin{aligned}
&|\det(S(M_{C}))|=|(\frac{1-\sqrt{3}\mathbf{i}}{2})^{b-a}-1|=|1-1|=0.\\
\end{aligned}
\end{equation*}
This completes the proof. $\Box$\\

Next we shall compute the determinant of the matrix $T(M_{C})$ for a mixed cycle $M_{C}$. For this purpose, we first give another division $\{\Psi_{1},\Psi_{2},\Psi_{3},\Psi_{4}\}$ for all mixed cycles. Let $a,b,c$ be still as described in the previous paragraph of Lemma 2.3. All mixed cycles with $a-b\equiv 1$ or $2$ or $4$ or $5(\mod 6)$ and odd number $c$ are divided into $\Psi_{1}$. All mixed cycles with $a-b\equiv 1$ or $2$ or $4$ or $5(\mod 6)$ and even number $c$ are divided into $\Psi_{2}$. All mixed cycles with $a-b\equiv 3$ or $0(\mod 6)$ and odd number $c$ are divided into $\Psi_{3}$. All mixed cycles with $a-b\equiv 3$ or $0(\mod 6)$ and even number $c$ are divided into $\Psi_{4}$. Moreover, if the cycle of a unicycle mixed graph belongs to $\Phi_{i}$ or $\Psi_{i}$ for some $i\in\{1,2,3,4\}$, then this unicycle mixed graph still belongs to it.\\
\\
\textbf{Lemma 2.4.} Let $M_{C}$ be a mixed cycle with $n$ vertices. Then
 \begin{equation*}
|\det(T(M_{C}))|=\left\{
\begin{array}{rcl}
&1,&\text{if $M_{C}$ belongs to $\Psi_{1}$};\\
&\sqrt{3},&\text{if $M_{C}$ belongs to $\Psi_{2}$};\\
&2,&\text{if $M_{C}$ belongs to $\Psi_{3}$};\\
&0,&\text{if $M_{C}$ belongs to $\Psi_{4}$}.\\
\end{array} \right.
\end{equation*}
\\
\textit{Proof.} Remark that the proof is similar to that of Lemma 2.3. We first expand the determinant of the matrix $T(M_{C})$ along the top row, that is,
\begin{equation*}
\begin{aligned}
|\det(T(M_{C}))|=&|(\frac{1-\sqrt{3}\mathbf{i}}{2})^{b}+(-1)^{n-1}(\frac{1-\sqrt{3}\mathbf{i}}{2})^{a}|\\
=&|(\frac{1-\sqrt{3}\mathbf{i}}{2})^{a}||(\frac{1-\sqrt{3}\mathbf{i}}{2})^{b-a}+(-1)^{n-1}|\\
=&|(\frac{1-\sqrt{3}\mathbf{i}}{2})^{b-a}+(-1)^{n-1}|.\\
\end{aligned}
\end{equation*}
First if $M_{C}$ belongs to $\Psi_{1}$, then $a-b\equiv 1$ or $2$ or $4$ or $5(\mod 6)$ and $c$ is odd. Thus, by a tedious calculation, one has
\begin{equation*}
\begin{aligned}
|\det(T(M_{C}))|=|(\frac{1-\sqrt{3}\mathbf{i}}{2})^{b-a}+(-1)^{n-1}|=1.\\
\end{aligned}
\end{equation*}
Similarly if $M_{C}$ belongs to $\Psi_{2}$, then $a-b\equiv 1$ or $2$ or $4$ or $5(\mod 6)$ and $c$ is even. So we get
\begin{equation*}
\begin{aligned}
|\det(T(M_{C}))|=|(\frac{1-\sqrt{3}\mathbf{i}}{2})^{b-a}+(-1)^{n-1}|=\sqrt{3}.\\
\end{aligned}
\end{equation*}
Now assume that $M_{C}$ belongs to $\Psi_{3}$. Then $a-b\equiv 3$ or $0(\mod 6)$ and $c$ is odd, which leads to
\begin{equation*}
\begin{aligned}
|\det(T(M_{C}))|=|(\frac{1-\sqrt{3}\mathbf{i}}{2})^{b-a}+(-1)^{n-1}|=2.\\
\end{aligned}
\end{equation*}
Finally if $M_{C}$ belongs to $\Psi_{4}$, then $a-b\equiv 3$ or $0(\mod 6)$ and $c$ is even. Therefore, we have
\begin{equation*}
\begin{aligned}
|\det(T(M_{C}))|=|(\frac{1-\sqrt{3}\mathbf{i}}{2})^{b-a}+(-1)^{n-1}|=0.\\
\end{aligned}
\end{equation*}
This completes the proof.   $\Box$\\

Let $M_{R}$ be a square substructure of a mixed graph. We call $M_{R}$ to be a substructure $\Rmnum{1}$ (abbreviate $S\Rmnum{1}$) if each of the components of $M_{R}$ is either a rootless tree, or a unicyclic graph that does not belong to $\Phi_{4}$. Similarly, we call $M_{R}$ to be a substructure $\Rmnum{2}$ (abbreviate $S\Rmnum{2}$) if each of the components of $M_{R}$ is either a rootless tree, or a unicyclic graph that does not belong to $\Psi_{4}$.\\
\\
\textbf{Lemma 2.5.} For a square substructure $M_{R}$ of a mixed graph, we have
\begin{equation*}
|\det (S(M_{R}))|=\left\{\begin{array}{rcl}
&\sqrt{3}^{\gamma_{1}}2^{\gamma_{2}},&\;\text{if $M_{R}$ is an $S\Rmnum{1}$};\\
&0,&\;\text{otherwise},\\
\end{array} \right.
\end{equation*}
where $\gamma_{1}$ and $\gamma_{2}$ are the number of unicyclic graphs of $\Phi_{2}$ and $\Phi_{3}$ in all the components of $S\Rmnum{1}$,
respectively.\\
\\
\textit{Proof.} First assume that $M_{R}$ is not an $S\Rmnum{1}$, then at least one component of $M_{R}$ is a unicyclic graph belonging to $\Phi_{4}$. According to Lemma 2.3, along with a simple analysis, we can get that $|\det (S(M_{R}))|=0$. Now assume that $M_{R}$ is an $S\Rmnum{1}$. Then there exists a permutation matrix $P$ such that $P^{-1}S(M_{R})P$ is a block diagonal matrix, in which each of the diagonal blocks corresponds to a component of $M_{R}$. It follows from the definition of the $S\Rmnum{1}$ that each of the components of $M_{R}$ is either a rootless tree, or a unicyclic graph that does not belong to $\Phi_{4}$. If $M_{R}$ has a component to be a rootless tree, then the absolute value of the determinant of the corresponding block is 1 from Lemma 2.1. If
$M_{R}$ has a component to be a unicyclic graph that does not belong to $\Phi_{4}$, then it clearly follows from Lemma 2.3 that the absolute value of the determinant of the corresponding block is 1, $\sqrt{3}$ and 2 according to the component belongs to $\Phi_{1}$, $\Phi_{2}$ and $\Phi_{3}$, respectively. Since the absolute value of $\det (S(M_{R}))$ is equal to the absolute value multiplication of the determinant of each diagonal block. Hence the required result holds. $\Box$\\
\\
\textbf{Lemma 2.6.} Let $M_{R}$ be a square substructure of a mixed graph. Then
\begin{equation*}
|\det (T(M_{R}))|=\left\{\begin{array}{rcl}
&\sqrt{3}^{\tau_{1}}2^{\tau_{2}},&\;\text{if $M_{R}$ is an $S\Rmnum{2}$};\\
&0,&\;\text{otherwise},\\
\end{array} \right.
\end{equation*}
where $\tau_{1}$ and $\tau_{2}$ are the number of unicyclic graphs of $\Psi_{2}$ and $\Psi_{3}$ in components of $S\Rmnum{2}$, respectively.\\
\\
\textit{Proof.} This proof is similar to that of Lemma 2.5, the detail is omitted. $\Box$\\

Given an $m$-by-$n$ matrix $A$, we use $\alpha$ to represent the subset of the set of its rows and $\beta$ to represent the subset of the set of its columns. Then $A[\alpha,\beta]$ denotes the submatrix of $A$ derived by taking the rows indexed by $\alpha$ and columns indexed by $\beta$. $A(\alpha,\beta)$ denotes the submatrix of $A$ obtained by deleting the rows corresponding to $\alpha$ and columns corresponding to $\beta$. Moreover, we simply write $A[\alpha]$ and $A(\alpha)$ in place of $A[\alpha,\alpha]$ and $A(\alpha,\alpha)$, respectively.\\
\\
\textbf{Theorem 2.7.} Let $M_{G}$ be a mixed graph with vertex set $V$ and edge set $E$, and $V_{1}$ be a vertex subset of $V$. Then
\[
\det (L[V_{1}])=\sum_{M_{R}}3^{\gamma_{1}}4^{\gamma_{2}},
\]
where the sum is over all the $M_{R}$ with the property $S\Rmnum{1}$ and $V(M_{R})=V_{1}$, and $\gamma_{1}$, $\gamma_{2}$ are the number of unicyclic graphs of $\Phi_{2}$ and $\Phi_{3}$ in components of $M_{R}$, respectively.\\
\\
\textit{Proof.} It is clear that $L[V_{1}]=S[V_{1},E]S[V_{1},E]^{\ast}$. Then the Cauchy-Binet Theorem implies that $\det L[V_{1}]$ is equal to the sum of absolute values of the squares of $\det S[V_{1},E_{1}]$, in which $E_{1}$ is an edge subset of $E$ satisfying $|V_{1}|=|E_{1}|$. So it follows from Lemma 2.5 that each substructure $M_{R}$ with the property $S\Rmnum{1}$ contributes $(\sqrt{3}^{\gamma_{1}}2^{\gamma_{2}})^{2}=3^{\gamma_{1}}4^{\gamma_{2}}$ to $\det (L[V_{1}])$, and the other substructures contribute 0 to $\det (L[V_{1}])$. This completes the proof. $\Box$\\
\\
\textbf{Theorem 2.8.} Let $M_{G}$ be a mixed graph with vertex set $V$ and edge set $E$, and $V_{1}$ be a vertex subset of $V$. Then
\[
\det (Q[V_{1}])=\sum_{M_{R}}3^{\tau_{1}}4^{\tau_{2}},
\]
where the sum is over all the substructure $M_{R}$ with the property $S\Rmnum{2}$ and $V(M_{R})=V_{1}$, and $\tau_{1}$ and $\tau_{2}$ are the number of unicyclic graphs of $\Psi_{2}$ and $\Psi_{3}$ in components of $M_{R}$, respectively.\\
\\
\textit{Proof}: Using a similar technique to Theorem 2.7, along with Lemma 2.6, we easily arrive at the desired result. $\Box$\\

Remarks that, from Theorem 2.7, we get directly that the Laplacian matrix of the second kind of a mixed graph is singular if and only if all its cycles (if possible) belong to $\Phi_{4}$. Similarly, Theorem 2.8 implies that the quasi-Laplacian matrix of the second kind of a mixed graph is singular if and only if all its cycles (if possible) belong to $\Psi_{4}$.

In what follows, we provide a structural characterization for the Laplacian matrix of the second kind being singular, which will be used to investigate non-principal minor version of matrix tree theorem for the Laplacian matrix of the second kind of a mixed graph. For this, we first introduce a special class of mixed graphs. Given a mixed graph $M_{G}$ with vertex set $V$ and edge set $E$, it is said to be $SP$ if $V$ can be divided into six (possibly empty) parts $V_{1}$, $V_{2}$, $V_{3}$, $V_{4}$, $V_{5}$ and $V_{6}$, such that all the undirected edges in $E$ are contained only in $V_{i}\;(i=1,2,3,4,5,6)$, while all the directed edges in $E$ are from $V_{i}$ to $V_{i+1}\;(i=1,2,3,4,5)$ and from $V_{6}$ to $V_{1}$, see Figure 3.
\begin{figure}[H]
\centering
\includegraphics[width=7.8cm,height=6.5cm]{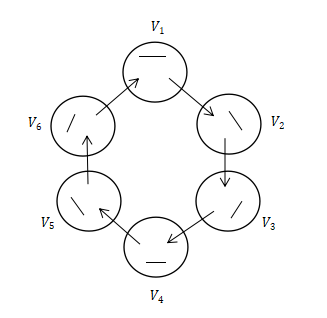}
\caption{The mixed graph with the structural property $SP$.}
\end{figure}

Let $W=v_{1}v_{2}\cdots v_{k-1}v_{k}$ be a mixed walk. We call a directed edge $e:v_{i}\rightarrow v_{j}$ is forward (resp. backward) in $W$ if $d_{v_{1}v_{j}}=d_{v_{1}v_{i}}+1$ (resp. $d_{v_{1}v_{j}}=d_{v_{1}v_{i}}-1$) in the underlying graph of $W$.\\
\\
\textbf{Lemma 2.9.} Let $M_{G}$ be an $SP$ mixed graph, and $u,v$ be two vertices in the same vertex subset. Then $a-b\equiv 0(\mod 6)$ for any mixed walk $W_{uv}$ connecting $u$ and $v$, where $a$ and $b$ represent the number of forward and backward edges in $W_{uv}$, respectively.\\
\\
\textit{Proof.} From the definition of $SP$ mixed graph, we can find that directed edges only exist between different vertex subsets, while undirected edges only exist in the same vertex subset. If all the vertices in $W_{uv}$ are in the same vertex set, then the result holds. Now suppose there is at least one vertex that is not in the same vertex subset as $u,v$. Then for each forward edge, there must be a backward edge corresponding to it, otherwise, the number of forward(backward) edges will be a miltiple of 6. Therefore, $a-b\equiv 0(\mod 6)$ will also be satisfied in the above situations. $\Box$\\
\\
\textbf{Theorem 2.10.} A mixed graph $M_{G}$ is $SP$ if and only if all the mixed cycles in $M_{G}$ belong to $\Phi_{4}$.\\
\\
\textit{Proof.} First suppose that all the mixed cycles in $M_{G}$ belong to $\Phi_{4}$. For any fixed vertex $v$, let $W_{uv}$ be any mixed walk in $M_{G}$ that connects any vertex $u\in V(M_G)$ and the fixed vertex $v$. Then $u$ can be divided into the following six vertex subsets in terms of any mixed walk $W_{uv}$.
\begin{equation*}
\begin{aligned}
V_{1}=&\{u\in V(M_{G}):a-b\equiv 1(\mod 6)\quad \text{for some $W_{uv}$}\},\\
V_{2}=&\{u\in V(M_{G}):a-b\equiv 2(\mod 6)\quad \text{for some $W_{uv}$}\},\\
V_{3}=&\{u\in V(M_{G}):a-b\equiv 3(\mod 6)\quad \text{for some $W_{uv}$}\},\\
V_{4}=&\{u\in V(M_{G}):a-b\equiv 4(\mod 6)\quad \text{for some $W_{uv}$}\},\\
V_{5}=&\{u\in V(M_{G}):a-b\equiv 5(\mod 6)\quad \text{for some $W_{uv}$}\},\\
V_{6}=&\{u\in V(M_{G}):a-b\equiv 0(\mod 6)\quad \text{for some $W_{uv}$}\}.\\
\end{aligned}
\end{equation*}
According to the vertex subsets divided above, we know that all vertices in $M_{G}$ are included in these six vertex subsets.

Now we first prove that $\pi: V_{1}\cup V_{2}\cup V_{3}\cup V_{4}\cup V_{5}\cup V_6$ constitutes a partition of vertex set $V(M_{G})$, equivalently, $V_i\cap V_j={\O}$ for $i\neq
j,1\leq i,j\leq 6$. Suppose that $u\in V_{1}\cap V_{2}$. Since $u\in V_{1}$, then there exists a walk $W_{1}$ with $a_{1}$ forward edges and $b_{1}$ backward edges so that $a_{1}-b_{1}\equiv 1(\mod 6)$. On the other hand, since $u\in V_{2}$, then there is a walk $W_{2}$ with $a_{2}$ forward edges and $b_{2}$ backward edges that satisfies $a_{2}-b_{2}\equiv 2(\mod 6)$. Obviously, the walks $W_{1}$ and $W_{2}$ form a mixed cycle $M_{C}$. If $W_{1}$ and $W_{2}$ have no common edge, then $M_{C}$ has $a=a_{1}+b_{2}$ forward edges and $b=b_{1}+a_{2}$ backward edges that satisfies  $a-b\equiv 5(\mod 6)$. So the mixed cycle $M_{C}$ belongs to $\Phi_{1}$, which contradicts that all the mixed cycles in $M_{G}$ belong to $\Phi_{4}$. If $W_{1}$ and $W_{2}$ have common edges, this result is also true by a simple analysis. Hence $V_{1}\cap V_{2}={\O}$. Similarly, we show that $V_i\cap V_j={\O}$ for $i\neq j,1\leq i,j\leq 6$.

In the following, we consider all the edges of $M_{G}$. Given a fixed edge $e_{uv}$, choose any vertex $w\in V(M_{G})$, let $W_{1}$ be a walk connecting $u$ and $w$ with $a_{1}$ forward edges and $b_{1}$ backward edges. Also let $W_{2}$ be a walk composed of $W_{1}$ and $e_{uv}$ with $a_{2}$ forward edges and $b_{2}$ backward edges. We divide it into two situations below.

\textit{Case 1:} $e_{uv}$ is an undirected edge.

Assume that $u\in V_{i_0}$ for some $i_0\in\{1,2,3,4,5,6\}$. Then for $W_{1}$ we have $a_{1}-b_{1}\equiv i_{0}(\mod 6)$. Bear in mind that $e_{uv}$ is an undirected edge, then $a_{2}-b_{2}\equiv i_{0}(\mod 6)$ for $W_{2}$. This indicates that $v\in V_{i_0}$. Hence, the undirected edge $e_{uv}$ is contained in vertex subset $V_{i_0}$.

\textit{Case 2:} $e_{uv}$ is a directed edge with a direction from $u$ to $v$.

Suppose that $u\in V_{1}$, then for $W_{1}$ we have $a_{1}-b_{1}\equiv 1(\mod 6)$. Since $a_{2}=a_{1}+1$ and $b_{2}=b_{1}$ for $W_{2}$, then $a_{2}-b_{2}\equiv 2(\mod 6)$. Thus we conclude that $v\in V_{2}$.
If $u\in V_{2}$, then $a_{1}-b_{1}\equiv 2(\mod 6)$ for $W_{1}$. For similar reasons, we obtain that $a_{2}-b_{2}\equiv 3(\mod 6)$ for $W_{2}$, which implies that $v\in V_{3}$.
If $u\in V_{3}$, then $a_{1}-b_{1}\equiv 3(\mod 6)$ for $W_{1}$, and $a_{2}-b_{2}\equiv 4(\mod 6)$ for $W_{2}$. So we get that $v\in V_{4}$.
If $u\in V_{4}$, then for $W_{1}$ we have $a_{1}-b_{1}\equiv 4(\mod 6)$, Thus $a_{2}-b_{2}\equiv 5(\mod 6)$ for $W_{2}$, which implies that $v\in V_{5}$.
If $u\in V_{5}$, then for $W_{1}$ we have $a_{1}-b_{1}\equiv 5(\mod 6)$, and $a_{2}-b_{2}\equiv 0(\mod 6)$ for $W_{2}$. So we conclude that $v\in V_{6}$.
If $u\in V_{6}$, then for $W_{1}$ we have $a_{1}-b_{1}\equiv 0(\mod 6)$, and $a_{2}-b_{2}\equiv 1(\mod 6)$ for $W_{2}$. So we conclude that $v\in V_{1}$.
Based on the above discussion, we get that the direction of any directed edge only occurs from $V_{i}$ to $V_{i+1}$ ($i=1,2,3,4,5$), or from $V_{6}$ to $V_{1}$.

Finally, by checking the definition of $SP$ mixed graph, we can conclude that $M_{G}$ is a mixed graph with property $SP$.

Conversely, assmue that the  mixed graph $M_{G}$ is $SP$. Let $M_{C}$ be any mixed cycle of $M_{G}$ and $u$ be any vertex on the cycle $M_{C}$. We may view $M_{C}$ as a mixed walk $W_{uu}$. Thus, $a-b\equiv 0(\mod 6)$ from Lemma 2.9. Hence, all the mixed cycles in $M_{G}$ belong to $\Phi_{4}$. This completes the proof. $\Box$\\

Combining Theorem 2.10 with Theorem 2.7, we get directly the following result, which gives a structural characterization of Laplacian matrix of the second kind being singular. Here we provide another alternative proof that is straight-forward and independent of the previous results.\\
\\
\textbf{Theorem 2.11.}  A mixed graph $M_{G}$ with $n$ vertices is $SP$ if and only if its Hermitian Laplacian matrix of the second kind $L(M_{G})$ is singular.\\
\\
\textit{Proof.} Suppose first that $L(M_{G})$ is singular. Since $L(M_{G})=S(M_{G})S(M_{G})^{\ast}$, then there is a vector $\xi=(\xi_{1},\xi_{2},\ldots,\xi_{n})^{\top}$ so that $\xi^{\ast}S(M_{G})=0$. Thus we obtain that $\overline{\xi_{i}}-\frac{1-\sqrt{3}\mathbf{i}}{2}\overline{\xi_{j}}=0$ for an edge $e:v_{i}\rightarrow v_{j}$, that is, $\xi_{j}=\frac{1-\sqrt{3}\mathbf{i}}{2}{\xi_{i}}$. Similarly, $\overline{\xi_{i}}-\overline{\xi_{j}}=0$ for an edge $e:v_{i}\leftrightarrow v_{j}$, which implies that $\xi_{j}=\xi_{i}$. Now we can divide the vertex set $V(M_{G})$ into the following six parts:
\begin{equation*}
\begin{aligned}
V_{1}=&\{v_{k}:\xi_{k}=\xi_{1}\},\\
V_{2}=&\{v_{k}:\xi_{k}=\frac{1-\sqrt{3}\mathbf{i}}{2}\xi_{1}\},\\
V_{3}=&\{v_{k}:\xi_{k}=-\frac{1+\sqrt{3}\mathbf{i}}{2}\xi_{1}\},\\
V_{4}=&\{v_{k}:\xi_{k}=-\xi_{1}\},\\
V_{5}=&\{v_{k}:\xi_{k}=-\frac{1-\sqrt{3}\mathbf{i}}{2}\xi_{1}\},\\
V_{6}=&\{v_{k}:\xi_{k}=\frac{1+\sqrt{3}\mathbf{i}}{2}\xi_{1}\}.\\
\end{aligned}
\end{equation*}
It is easy to see that $\pi:V_{1}\cup V_{2}\cup V_{3}\cup V_{4}\cup V_{5}\cup V_{6}$ is a partition of $V(M_{G})$ and satisfies the definition of $SP$ mixed graph.

Conversely, assume that $M_{G}$ is an $SP$ mixed graph. As a matter of convenience, let $n_i=|V_{i}|$ for $i=1,2,\dots,6$. Let $\xi$ be an $n$-dimensional column vector such that
\begin{equation*}
\begin{aligned}
\xi=\left(\mathbf{1}_{n_1}^{\top},\frac{1-\sqrt{3}\mathbf{i}}{2}\mathbf{1}_{n_2}^{\top},-\frac{1+\sqrt{3}\mathbf{i}}{2}\mathbf{1}_{n_3}^{\top},-\mathbf{1}_{n_4}^{\top},-\frac{1-\sqrt{3}\mathbf{i}}{2}\mathbf{1}_{n_5}^{\top},\frac{1+\sqrt{3}\mathbf{i}}{2}\mathbf{1}_{n_6}^{\top}\right)^{\top},
\end{aligned}
\end{equation*}
where $\mathbf{1}_{n_1}^{\top}$ denotes the transpose of the column vector of order $n_1$ with all entries 1. It is easy to verify that $\xi^{\ast}S(M_{G})=0$. Therefore, $\xi$ is an eigenvector for $L(M_{G})$ corresponding to the eigenvalue $0$, which implies $L(M_{G})$ is singular. $\Box$

\section{Non-principal minors}

The following definitions and notations are introduced in the reference \cite{Bapat}. Let $M_{G}$ be a mixed graph with vertex set $V$ and edge set $E$. If $V'\subseteq V$, $E'\subseteq E$, then the substructure of $M_{G}$ consisting of the vertices in $V'$ and edges in $E'$ is denoted by $\mathbb{S}(V',E')$, whose incidence matrix corresponds to the submatrix $S[V',E']$ of $S(M_G)$. We say  $\mathbb{S}(V',E')$ to be non-singular relative to $V'$ if the matrix $S[V',E']$ is non-singular, where $|V'|=|E'|$. Assume that $V_1, V_2\subseteq V$, $E_1\subseteq E$ with $|V_{1}|=|V_{2}|=|E_{1}|$. We call $\mathbb{S}(V_{1}\cup V_{2},E_{1})$ is non-singular relative to $V_{1}$ and $V_{2}$ if $S[V_{1},E_{1}]$ and $S[V_{2},E_{1}]$ are both non-singular. Under the circumstances, we also call $\mathbb{S}(V_{1}\cup V_{2},E_{1})$ a generalized matching between $V_{1}$ and $V_{2}$. Moreover, if $L[V_{1},V_{2}]$ is a off-diagonal square submatrix of $L(M_{G})$, then $L[V_{1},V_{2}]=S[V_{1},E_{1}]S[V_{2},E_{1}]^{\ast}$.

Similarly, for $V'\subseteq V$, $E'\subseteq E$, we use $\mathbb{T}(V',E')$ to denote the substructure of $M_{G}$ corresponding to the submatrix $T[V',E']$ of $T(M_G)$. Thus, for $V_1, V_2\subseteq V$, $E_1\subseteq E$ with $|V_{1}|=|V_{2}|=|E_{1}|$, a non-singular substructure $\mathbb{T}(V_{1}\cup V_{2},E_{1})$ is said to be a quasi-generalized matching between $V_{1}$ and $V_{2}$ if $T[V_{1},E_{1}]$ and $T[V_{2},E_{1}]$ are both non-singular. In the same way, if $Q[V_{1},V_{2}]$ is a off-diagonal square submatrix of $Q(M_{G})$, then  $Q[V_{1},V_{2}]=T[V_{1},E_{1}]T[V_{2},E_{1}]^{\ast}$.

We first recall the following two lemmas, which are due to Bapat, Grossman and Kulkarni \cite{Bapat}.\\
\\
\textbf{Lemma 3.1.}\cite{Bapat} Let $M_{G}$ be a mixed graph with vertex set $V$ and edge set $E$, and $V_{1},V_{2}\subseteq V$, $E_{1}\subseteq E$ with $|V_{1}|=|V_{2}|=|E_{1}|$. Then $\mathbb{S}(V_{1}\cup V_{2},E_{1})$ is non-singular relative to $V_{1}$ and $V_{2}$ if and only if each component is either a non-singular substructure of $\mathbb{S}(V_{1}\cap V_{2},E_{1})$ or a tree with exactly one vertex in each of $V_{1}\setminus V_{2}$ and $V_{2}\setminus V_{1}$.\\
\\
\textbf{Lemma 3.2.}\cite{Bapat} Let $M_{G}$ be a mixed graph with vertex set $V$ and edge set $E$, and $V_{1},V_{2}\subseteq V$, $E_{1}\subseteq E$ with $|V_{1}|=|V_{2}|=|E_{1}|$. Then $\mathbb{T}(V_{1}\cup V_{2},E_{1})$ is non-singular relative to $V_{1}$ and $V_{2}$ if and only if each component is either a non-singular substructure of $\mathbb{T}(V_{1}\cap V_{2},E_{1})$ or a tree with exactly one vertex in each of $V_{1}\setminus V_{2}$ and $V_{2}\setminus V_{1}$.\\
\\
\textbf{Lemma 3.3.} Let $M_{T}$ be a mixed tree in a generalized matching $\mathbb{S}(V_{1}\cup V_{2},E_{1})$ between $V_{1}$ and $V_{2}$, with edge vertex $E_{1}$. Then $M_{T}$ contributes
\[
(-\frac{1+\sqrt{3}\mathbf{i}}{2})^{b-a}\cdot (-1)^{c}
\]
to $\det L[V_{1},V_{2}]$, where $a$ and $b$ represent the respective number of directed edges away from the vertices of $V_{1}\setminus V_{2}$ and $V_{2}\setminus V_{1}$ in the path connecting the two vertices in $M_{T}$, and $c$ represents the number of undirected edges in this path.\\
\\
\textit{Proof.} Since $M_{T}$ is a mixed tree in a generalized matching $\mathbb{S}(V_{1}\cup V_{2},E_{1})$ between $V_{1}$ and $V_{2}$, then there is exactly one vertex in each of $V_{1}\setminus V_{2}$ and $V_{2}\setminus V_{1}$ by Lemma 3.1. Let $u\in V_{1}\setminus V_{2}$ and $v\in V_{2}\setminus V_{1}$ be the two vertices in $M_{T}$. Then there exists a rootless tree $M_{T}\setminus \{v\}$ that is a component of an $S\Rmnum{1}$ corresponding to $V_{1}$ with $|V_{1}|=|E_{1}|$. Likewise, there exists a rootless tree $M_{T}\setminus \{u\}$ that is a component of an $S\Rmnum{1}$ corresponding to $V_{2}$ with $|V_{2}|=|E_{1}|$. By Lemma 2.1, we can get that $M_{T}\setminus \{v\}$ and $M_{T}\setminus \{u\}$ contribute $(-\frac{1-\sqrt{3}\mathbf{i}}{2})^{a_{1}}\cdot (-1)^{\beta_{1}'}$ and $(-\frac{1-\sqrt{3}\mathbf{i}}{2})^{b_{1}}\cdot (-1)^{\beta_{2}'}$ to $\det S[V_{1},E_{1}]$ and $\det S[V_{2},E_{1}]$, respectively, where $a_{1}$ and $b_{1}$ represent the respective number of old directed edges away from $v$ and $u$, $\beta_{1}'$ and $\beta_{2}'$ are the respective number of new directed edges away from $v$ and $u$ after we assign any direction to each undirected edge in $M_{R}$. Thus, by $L[V_{1},V_{2}]=S[V_{1},E_{1}]S[V_{2},E_{1}]^{\ast}$, we find that $M_{T}\setminus \{v\}$ and $M_{T}\setminus \{u\}$ contribute
\[
(-\frac{1-\sqrt{3}\mathbf{i}}{2})^{a_{1}}\cdot \overline{(-\frac{1-\sqrt{3}\mathbf{i}}{2})^{b_{1}}}\cdot (-1)^{\beta_{1}'+\beta_{2}'}=(-\frac{1+\sqrt{3}\mathbf{i}}{2})^{b_{1}-a_{1}}\cdot (-1)^{\beta_{1}'+\beta_{2}'}
\]
to $\det L[V_{1},V_{2}]$. But for any edge that is not on the path connecting $u$ and $v$, it is always away from or towards $u$ and $v$ at the same time. In either case, its contribution to $\det L[V_{1},V_{2}]$ is $1$, it is to say, we can ignore all the edges that are not on the path connecting $u$ and $v$. Therefore, the mixed tree $M_{T}$ contributes
\[
(-\frac{1+\sqrt{3}\mathbf{i}}{2})^{b-a}\cdot (-1)^{\beta_{1}+\beta_{2}}=(-\frac{1+\sqrt{3}\mathbf{i}}{2})^{b-a}\cdot (-1)^{c}
\]
to $\det L[V_{1},V_{2}]$, where $a$ and $b$ represent the respective number of old directed edges away from $v$ and $u$ in the path connecting $v$ and $u$ in $M_{T}$, $\beta_{1}$ and $\beta_{2}$ are the respective number of new directed edges away from $v$ and $u$ after we assign any direction to all the undirected edges in this path, and $c$ is the number of undirected edges in this path. This completes the proof. $\Box$\\
\\
\textbf{Lemma 3.4.} Let $M_{T}$ be a mixed tree in a quasi-generalized matching $\mathbb{T}(V_{1}\cup V_{2},E_{1})$ between $V_{1}$ and $V_{2}$. Then  $M_{T}$ contributes
\[
(\frac{1+\sqrt{3}\mathbf{i}}{2})^{b-a}
\]
to $\det Q[V_{1},V_{2}]$, where $a$ and $b$ represent the number of directed edges away from the two vertices of $V_{1}\setminus V_{2}$ and $V_{2}\setminus V_{1}$ in the path connecting them in $M_{T}$, respectively.\\
\\
\textit{Proof.} Using exactly the same method as Lemma 3.3, along with Lemmas 2.2 and 3.2, we can obtain the desired result. Some details of the proof are omitted. $\Box$\\

In the following, we shall provide the main results of this section, which give the explicit expressions of non-principal minors for Hermitian (quasi-)Laplacian matrix of the second kind of a mixed graph $M_{G}$.\\
\\
\textbf{Theorem 3.5.} Let $M_{G}$ be a mixed graph with vertex set $V$ and edge set $E$, and $V_{1},V_{2}\subseteq V$ with $|V_{1}|=|V_{2}|$. Then
\[
\det L[V_{1},V_{2}]=\sum_{\Gamma}(-\frac{1+\sqrt{3}\mathbf{i}}{2})^{\sum_{M_{T}}(b-a)}\cdot (-1)^{\sum_{M_{T}}c}\cdot 3^{\gamma_{1}}\cdot 4^{\gamma_{2}},
\]
where the first sum is taken over all generalized matching $\Gamma$ between $V_{1}$ and $V_{2}$ and the second sum is taken over all trees $M_{T}$ in $\Gamma$, $a$ and $b$ represent the respective number of directed edges away from the two vertices of $V_{1}\setminus V_{2}$ and $V_{2}\setminus V_{1}$ in the path connecting them in $M_{T}$, $c$ is the number of undirected edges in the path, $\gamma_{1}$ and $\gamma_{2}$ are the number of unicyclic graphs of $\Phi_{2}$ and $\Phi_{3}$ in components of $\Gamma$, respectively.\\
\\
\textit{Proof.} First Lemma 3.1 implies that there is either a non-singular substructure of $\mathbb{S}(V_{1}\cap V_{2},E_{1})$ or a tree with exactly one vertex in each of $V_{1}\setminus V_{2}$ and $V_{2}\setminus V_{1}$ for each component of a generalized matching $\Gamma$ between $V_{1}$ and $V_{2}$. So it follows from Theorem 2.7 and Lemma 3.3 that the required result holds. $\Box$\\
\\
\textbf{Theorem 3.6.} Let $M_{G}$ be a mixed graph with vertex set $V$ and edge set $E$, and $V_{1},V_{2}\subseteq V$ with $|V_{1}|=|V_{2}|$. Then
\[
\det Q[V_{1},V_{2}]=\sum_{\Gamma}(\frac{1+\sqrt{3}\mathbf{i}}{2})^{\sum_{M_{T}}(b-a)}\cdot 3^{\tau_{1}}\cdot 4^{\tau_{2}},
\]
where the first sum is taken over all quasi-generalized matching $\Gamma$ between $V_{1}$ and $V_{2}$ and the second sum is taken over all trees $M_{T}$ in $\Gamma$, $a$ and $b$ represent the number of directed edges away from the two vertices of $V_{1}\setminus V_{2}$ and $V_{2}\setminus V_{1}$ in the path connecting them in $M_{T}$, $\tau_{1}$ and $\tau_{2}$ are the number of unicyclic graphs of $\Psi_{2}$ and $\Psi_{3}$ in components of $\Gamma$, respectively.\\
\\
\textit{Proof.} Using Lemma 3.2, along with Theorem 2.8 and Lemma 3.4, we easily obtain the required result in a similar way to Theorem 3.5. $\Box$\\

Finally, we give the generalized matrix tree theorem for Hermitian (quasi-)Laplacian matrix of the second kind of mixed graphs before closing this paper.\\
\\
\textbf{Theorem 3.7.} Let $M_{G}$ be an $SP$ mixed graph. Then the absolute value of all the cofactors of the second kind of Hermitian Laplacian matrix $L(M_{G})$ are equal, and their common absolute value is the number of spanning trees of the underlying graph $G$.\\
\\
\textit{Proof.} The result follows from Theorems 2.7 and 3.5. $\Box$\\
\\
\textbf{Theorem 3.8.} Let $M_{G}$ be a mixed graph with all cycles belonging to $\Psi_{4}$. Then the absolute value of all the cofactors of the second kind of Hermitian quasi-Laplacian matrix $Q(M_{G})$ are equal, and their common absolute value is the number of spanning trees of the underlying graph $G$.\\
\\
\textit{Proof.} The result follows from Theorems 2.8 and 3.6. $\Box$\\
\\
\textbf{Remarks.} Recall that, for any unoriented graph $G$, we can calculate the number of its spanning trees in terms of the cofactors of its Laplacian matrix. But in general, we can't compute the number of spanning trees in terms of the cofactors of its quasi-Laplacian matrix. Here, we give two sufficient conditions that, for some mixed graphs, the absolute values of all the cofactors of Hermitian (quasi-)Laplacian matrix of the second kind are equal to the number of spanning trees of the underlying graph $G$. It is interesting to note that, for some special mixed graphs, the number of spanning trees of their underlying graphs can be calculated according to the cofactors of Hermitian (quasi-)Laplacian matrix of the second kind, which is in general different from that of  unoriented graphs. Now let's go back to Examples 1.1 and 1.2 in this article, in Example 1.1, let $v_{1}$ and $v_{4}$ belong to vertex subset $V_{1}$, $v_{2}$ and $v_{3}$ belong to vertex subset $V_{2}$ and otherwise, empty for mixed graph $M_{G}'$. It is easy to check that $M_{G}'$ is an $SP$ mixed graph. Applying Theorem 3.7, we can get the number of spanning trees of the underlying graph $G$ by calculating the absolute value of any cofactor of Hermitian Laplacian matrix of the second kind $L(M_{G}')$. Similarly, by a simple observation, we can find that all cycles of $M_{G}'$ in Example 1.2 belong to $\Psi_{4}$. Applying Theorem 3.8, we can get the number of spanning trees of its underlying graph $G$ by calculating the absolute value of any cofactor of Hermitian quasi-Laplacian matrix of the second kind $Q(M_{G}')$. These are just as we expected.

\end{document}